\documentclass[11pt]{amsart}

\usepackage{amssymb,amsmath,amsxtra,amscd}
\usepackage[mathscr]{eucal} 

\def\QQ         {{\bf Q}}

\def\ii         {{\rm i}}
\def\ee         {{\rm e}}

\def\dim        {{\rm dim}}

\newcommand{\nonu}{\nonumber \\[2mm]}

\theoremstyle{plain}
\numberwithin{equation}{section}
\newtheorem{theorem}{Theorem}

\newtheorem{corollary}{Corollary}
\newtheorem{lemma}{Lemma}
\numberwithin{lemma}{section}
\numberwithin{theorem}{section}
\numberwithin{corollary}{section}
\numberwithin{proposition}{section}

\theoremstyle{definition}
\newtheorem{definition}{Definition}
\numberwithin{definition}{section}

\theoremstyle{remark}

\numberwithin{example}{section}

\newcommand{\nc}{\newcommand}
\nc{\mc}{\mathcal}
\nc{\MSV}[1]{\ensuremath{ \Omega^{\on{ch}}_{#1 } }}
\nc{\MSVg}[1]{\ensuremath{ \Omega^{\on{ch},g}_{#1 } }}
\nc{\on}{\operatorname}
\nc{\BunG}{\on{BunG}}
\nc{\BunGH}{\on{\BunG^{H}}}
\nc{\Mg}{\on{M}^{g}}
\nc{\M}{\on{M}}

\nc{\Cl}{Cl}
\nc{\drva}{\widehat{\Omega}}
\nc{\spec}{\on{Spec}}
\nc{\AutO}{\on{Aut} \mc{O}}
\nc{\vac}{|0\rangle}
\nc{\Z}{\mathbb{Z}}
\nc{\zf}[1]{z^{\frac{1}{#1}}}
\nc{\wf}[1]{w^{\frac{1}{#1}}}
\nc{\Mt}{M^{\sigma}}
\renewcommand{\exp}{\on{exp}}
\nc{\gr}{\on{gr}}

\newcommand{\is}{& \!\! = \!\! &}

\nc{\al}{\alpha}
\nc{\ol}{\overline}

\begin{document}

\title{Discrete torsion, orbifold elliptic genera, and the chiral de Rham complex}

\author{Anatoly Libgober} \thanks{}
\address{Department of Mathematics  
         University of Illinois at Chicago, Chicago, IL 60637}
\email{libgober@math.uic.edu}
\author{Matthew Szczesny}
\address{Department of Mathematics  
         University of Pennsylvania, Philadelphia, PA 19104 }
\email{szczesny@math.upenn.edu}

\date{December 2004}

\begin{abstract}

Given a compact complex algebraic variety with an effective action of a finite
group $G$, and a class $\alpha \in H^2(G,U(1))$, we introduce an
orbifold elliptic genus with discrete torsion $\alpha$, denoted
$Ell^{\alpha}_{orb}(X,G, q, y)$. We give an interpretation of this
genus in terms of the chiral de Rham complex attached to the orbifold
$[X/G]$. If $X$ is Calabi-Yau and $G$ preserves the
volume form, $Ell^{\alpha}_{orb}(X,G, q, y)$ is a weak Jacobi form. We also
obtain a formula for the generating function of the elliptic genera of
symmetric products with discrete torsion.

\end{abstract}

\maketitle 

\section{Introduction}

The two-variable elliptic genus (see for example \cite{Kri}) of a
compact complex manifold $X$ is a generating function
\begin{equation} \label{ellgen}
Ell(X,y,q) = \sum_{m,l} c(m,l) q^{m} y^{l}
\end{equation}
which captures important topological information about $X$. For
appropriate values of $y$ and $q$, $Ell(X,y,q)$ specializes to the L,
$\widehat{A}$, and $\chi_y$ genera respectively.  Mathematically, the
elliptic genus can be defined as follows. For a holomorphic vector
bundle $V$ on $X$ and a formal variable $t$, let
\[
\on{Sym}_t V = 1 + t V + t^2 \on{Sym}^{2} V + t^3 \on{Sym}^{3} V +
\cdot \; \; \in K_{0} (X)[[t]]
\]
and
\[
\Lambda_t V = 1 + t V + t^2 \Lambda^{2} V + t^3 \Lambda^{3} V + \cdot
\; \; \in K_{0} (X)[[t]]
\]
Let $T_X, T^{*}_X$ denote the holomorphic tangent and cotangent bundles respectively, and 
\begin{equation}
{\mc Ell}(X,q,y) = y^{-{\dim X \over 2}}
 \otimes_{n \ge 1} (\Lambda_{-yq^{n-1}}
T^{*}_X \otimes \Lambda_{-y^{-1}q^n} T_{X} \otimes S_{q^n} T^{*}_X \otimes
S_{q^n} T_X) 
\label{bundle}
\end{equation}
viewed as an element of $K_{0}(X)[[q]][[y^{\pm \frac{1}{2}} ]]$. Then
\[
Ell(X,q,y) = \chi({\mc Ell}(X,q,y)).
\]
In physics, $Ell(X,q,y)$ is part of the partition function of a
two-dimensional conformal field theory with target $X$. 

In this paper, we will say that $X$ is Calabi-Yau if $K_X$ is trivial - this is of course weaker than the usual mathematical Calabi-Yau condition, but  agrees with the Physics terminology.
When $X$ is Calabi-Yau,
$Ell(X,q,y)$ has nice modular properties. Let $\mathbb{H}$ denote the
upper half plane. A \emph{weak Jacobi form} of weight $k \in
\mathbb{Z}$ and index $r \in \frac{1}{2} \mathbb{Z}$ is a holomorphic
function on $\mathbb{H} \times \mathbb{C}$ satisfying the
transformation property
\[
\phi({{a\tau +b} \over {c \tau +d}},{z \over {c \tau +d}})=
 (c \tau +d)^k e^{2 \pi i{{r cz^2} \over {c \tau +d}}  }\phi (\tau, z), \; \; \; \; \left( \begin{matrix} a & b \\ c & d \end{matrix} \right) \in SL(2,\mathbb{Z})
\]
\[
\phi(\tau, z+m\tau+n)=
(-1)^{2r(m+n)} e^{-2 \pi i r (m^2 \tau +2 m z)} \phi (\tau, z), \; \; \; \; (m,n) \in \mathbb{Z}^2
\]
that in addition has a Fourier expansion $\sum_{l,m}c_{m,l}y^lq^m$
with nonnegative $m$ (see \cite{EZ}), where $q=e^{2 \pi i \tau}, y = e^{2 \pi i z}$. 
It is shown in \cite{BL, Gri} that if $X$ is
Calabi-Yau then $Ell(X,q,y)$ is a weak Jacobi
form of weight $0$ and index $\on{dim}(X)/2$.

In \cite{BL1}, the authors introduced a notion of orbifold elliptic
genus $Ell(X,G,q,y)$ attached to the global quotient orbifold $[X/G]$,
where $X$ is a smooth compact complex manifold, and $G \in Aut(X)$ is
a finite group (see also \cite{DLiM} for general reduced orbifolds). 
A mathematical definition of this genus is given in
section \ref{OrbGen}. Roughly, $Ell(X,G,q,y)$ is obtained by adding the
contributions of Euler characteristics of bundles analogous to
\ref{bundle} over the various fixed-point sets $X^g, g \in G$ of the
$G$--action on $X$. Furthermore, for each $g \in G$, the contribution
takes into account the eigenvalues of $g$ on $TX \vert_{X^{g}}$.  It
is also shown in \cite{BL1} that if $X$ is Calabi-Yau, and $G$
preserves the volume form, then $Ell(X,G,q,y)$ is a weak Jacobi form
of weight $0$ and index $\on{dim}(X)/2$.

It was observed in  \cite{BL} that the ordinary elliptic genus
$Ell(X,y,q)$ can be interpreted in terms of an object
called the { \em chiral de Rham complex}, denoted $\Omega^{ch}_X$. The
latter is a sheaf of vertex superalgebras attached to any smooth
complex manifold $X$, introduced in \cite{MSV} (for a brief discussion
of vertex algebras, see section \ref{VA}).  The sheaf $\Omega^{ch}_X$
has an increasing filtration $F^{\bullet} \Omega^{ch}_X$, and also a
compatible bigrading by two operators $L_0, J_0$. One can therefore
describe the associated graded sheaf $gr_F \Omega^{ch}_X$ in terms of the
bigraded components. One finds that $\on{Supertrace}(q^{L_0} y^{J_0},
\Omega^{ch}_X)$ is the sheaf \ref{bundle} above. It follows from this
observation that
\[
Ell(X,y,q) = \on{Supertrace} (q^{L_0} y^{J_0 - \on{dim} (X)/2}, H^{*}(X,\Omega^{ch}_X))
\] 
It is believed that $H^{*}(X, \Omega^{ch}_X)$ captures some of the
information of the two-dimensional quantum field theory on $X$
mentioned above. $\Omega^{ch}_X$ also carries a differential $Q_{BRST}
\in \on{End} (\Omega^{ch}_X)$, $Q^{2}_{BRST} = 0$ (which is why it is
called a complex). The "de Rham" part of the name comes from the fact
that the complex $(\Omega^{ch}_X, Q_{BRST})$ is quasiisomorphic to the
holomorphic de Rham complex $(\Omega_{dR}, \partial)$.

In \cite{FS}, the construction of $\Omega^{ch}_X$ was extended to
orbifolds (another construction of the chiral de
Rham complex for orbifolds was obtained independently by
A. Vaintrob). For each $g \in G$, one introduces
sheaves $\Omega^{ch,g}_X$ supported on $X^{g}$, which 
are modules over $\Omega^{ch}_X$. Each $\Omega^{ch,g}_X$ carries
a canonical $C(g)$--equivariant structure, where $C(g)$ denotes the
centralizer of $g$ in $G$. The sheaves $\Omega^{ch,g}_X$ allow one to
interpret some of the "stringy" geometric invariants of the orbifold
$[X/G]$. In particular, it is shown in \cite{FS} that
\begin{equation} \label{ellCDR}
Ell(X,G,q,y) = \on{Supertrace} (q^{L_0} y^{J_0 - \on{dim} (X)/2}, \bigoplus_{[g]} H^{*}
(X, \Omega^{ch,g}_X)^{C(g)} )
\end{equation}
and
\begin{equation} \label{CR}
\bigoplus_{[g]} \mathbb{H}^{*} (\Omega^{ch,g}_X, Q_BRST)^{C(g)} \cong \bigoplus_{[g]} H^{*}_{dR} (X^{g}/C(g), \mathbb{C})
\end{equation}
where $[g]$ denotes a set of representatives for the conjugacy classes in $G$.
The object on the right in \ref{CR}, with an appropriate grading and
ring structure is called the Chen-Ruan cohomology of $[X/G]$ (see
\cite{CR}). The isomorphism \ref{CR} is an isomorphism of graded
vector-spaces.

We now come to \emph{discrete torsion}. This term was introduced in
the physics literature to refer to the discovery (see \cite{Va1, Va2})
that an orbifold quantum field theory on $[X/G]$ could be "twisted" by
a cocycle $\alpha \in H^{2}(G,U(1))$. In terms of Physics, the path
integral can be written as a sum of contributions from sectors
parametrized by commuting pairs of elements $(g,h)$ in $G \times G$,
and the contribution from the $(g,h)$--sector is multiplied by the
phase $\delta(g,h) = \alpha(g,h)/ \alpha(h,g)$. This modification
produces a consistent physical theory and leads to $\alpha$--twisted
invariants of the orbifold $[X/G]$. For a mathematical discussion of
various aspects of discrete torsion see \cite{AF, R, Ka}.  In this
paper, we build on the results in \cite{BL1, FS} to give a
mathematical treatment of orbifold elliptic genera with discrete
torsion $Ell^{\alpha} (X,G,q,y)$. One way to define this object along
the lines of \cite{BL1} is to multiply each contribution in the sum
\ref{commute} by the appropriate phase $\delta(g,h)$. From the point
of view of the chiral de Rham complex, this definition can be cast in
a manner somewhat closer to the original physics approach as
follows. Recall that in \ref{ellCDR} one uses the $C(g)$--equivariant
structure on $\Omega^{ch,g}_X$ to project on the $C(g)$--invariant
part of $H^{*}(X, \Omega^{ch,g}_X)$. A cocycle $\alpha \in
H^{2}(G,U(1))$ induces characters $\alpha_g : C(g) \mapsto U(1)$ by
$\alpha_g (h) = \delta(g,h)$. The character $\alpha_g$ allows us to
twist the $C(g)$ equivariant structure, and taking invariants with
respect to this twisted structure projects on a different subspace. We
can therefore define
\begin{equation}
Ell^{\alpha} (X,G,q,y) = \on{Supertrace} (q^{L_0} y^{J_0 - \on{dim} (X)/2},
\bigoplus_{[g]} H^{*} (X, \Omega^{ch,g}_X)^{C(g)_{\alpha}} )
\end{equation}
where $C(g)_{\alpha}$ indicates that the twisted action is being
used. We show that if $X$ is Calabi-Yau, $Ell^{\alpha}
(X,G,q,y)$ is a Jacobi form of weight $0$ and index $\on{dim}(X)/2$. We
also show that there is an isomorphism of graded vector spaces
\[
\bigoplus_{[g]} \mathbb{H}^{*} (\Omega^{ch,g}_X, Q_{BRST})^{C(g)_{\alpha}} \cong
\bigoplus_{[g]} H^{*}_{dR} (X^{g}/C(g), \mc{L}^{\alpha}_g)
\]
where the object on the right denotes the $\alpha$--twisted Chen-Ruan
cohomology  of $[X/G]$ valued in the collection of local systems $\mc{L}^{\alpha}_g$, introduced in \cite{R}.

An important example of discrete torsion arises in the case of
symmetric products (see \cite{Di}). $S_N$, the symmetric group on $N$ letters, acts on
the hyperplane in $\mathbb{R}^{N}$ given by the equation $x_1 + \cdots
+ x_N = 0$. This yields an embedding $S_N \hookrightarrow
O(N-1)$. Pulling back the double cover $Pin(N-1) \mapsto O(N-1)$ yields a
central extension of $S_N$ by $\mathbb{Z} / 2 \mathbb{Z}$, which we
denote $\widehat{S}_N$ - i.e.
\begin{equation} \label{spincover}
1 \mapsto \mathbb{Z} / 2 \mathbb{Z} \mapsto \widehat{S}_N \mapsto S_N \mapsto 1
\end{equation}
The extension \ref{spincover} is non-split for $N \geq 4$, and
therefore yields a non-zero class $\alpha \in H^{2}(S_N,\mathbb{Z} /
2 \mathbb{Z} )$, which via the inclusion $\mathbb{Z}/2 \mathbb{Z}
\hookrightarrow U(1)$ can be pushed into $H^{2}(S_N,U(1))$. 

The orbifold elliptic genera of symmetric products can be arranged
into remarkable generating functions. It was proved in \cite{BL1}
following a physics derivation in \cite{DMVV} that
\[
\sum_{N \in \mathbb{Z_{+}}} p^{N} Ell_{orb} (X^{N}, S_N, y, q) =
 \prod_{n,m,\ell \geq 0} (1-p^nq^my^\ell)^{-c(nm,\ell)}.
\]
where the $c(m,l)$'s are as in \ref{ellgen}. In section \ref{symtors}
we obtain a generalization of this formula with discrete torsion given
by $\alpha$ above, which was originally obtained by Dijkgraaf
(\cite{Di}) in the physics literature.

{\bf Acknowledgements:} This project originally began with Lev
Borisov. We would like to thank him for many valuable
conversations. During the course of this work the second author was
supported by NSF grant $\on{DMS}-0401619$.

\section{Orbifold elliptic genera} 

\subsection{The orbifold elliptic genus} \label{OrbGen}

Let $X$ be a complex manifold on which a finite group $G$
acts effectively via holomorphic transformations. 
Let $X^h$ will be the fixed point set of $h \in G$
and  $X^{g,h}=X^g \cap X^h (g,h \in G)$. Let 
\begin{equation}
TX \vert_{X^h}= \oplus_{\lambda(h) \in \QQ \cap [0,1)} V_{\lambda}.
\label{directsum}
\end{equation}
where the bundle $V_{\lambda}$ on $X^h$ is determined by the 
requirement that $h$ acts on $V_{\lambda}$ via multiplication by 
$e^{2 \pi i \lambda(h)}$. For a connected component of $X^h$
(which by abuse of notation we also will denote $X^h$), the fermionic 
shift is defined as $F(h,X^h \subseteq X)=\sum_{\lambda} \lambda(h)$
(cf. \cite{Z}, \cite{BD}). Let us consider the 
bundle:
$$
V_{h,X^h\subseteq X}:= \otimes_{k\geq 1} 
\Bigr[
(\Lambda^\bullet_{yq^{k-1}} V_0^*)\otimes
(\Lambda^\bullet_{y^{-1} q^k} V_0)\otimes
(Sym^\bullet_{q^{k}} V_0^*)\otimes
(Sym^\bullet_{q^{k}} V_0 )\otimes
$$
\begin{equation}\otimes
\bigl[
\otimes_{\lambda\neq 0} 
(\Lambda^\bullet_{ yq^{k-1+\lambda(h)}} V_\lambda^* )\otimes
(\Lambda^\bullet_{ y^{-1}q^{k-\lambda(h)}} V_\lambda)\otimes
(Sym^\bullet_{q^{k-1+\lambda(h)}} V_\lambda^*)\otimes
(Sym^\bullet_{q^{k-\lambda(h)}} V_\lambda)
\bigr]
\Bigl]
\end{equation}
\begin{definition} {\rm The orbifold elliptic genus of a $G$-manifold $X$
is the function on $H \times {\bf C}$ given by:
$$Ell_{orb} (X,G,q,y) := 
y^{-\dim  X /2} \sum_{[g],X^g}   y^{F(g,X^g\subseteq X)} \frac 1 {|C(g)|}
\sum_{h\in C(h)} L(h, V_{g,X^g\subseteq X})
$$ where the summation in the first sum is over all conjugacy classes
in $G$ and connected components $X^g$ of an element $g \in [g]$,
$C(g)$ is the centralizer of $g \in G$ and
\[ 
L(h,V_{g,X^g\subseteq X})=
\sum_i (-1)^i{\rm tr}(h, H^i(V_{g,X^g\subseteq X}))
\]
 is the holomorphic
Lefschetz number.}
\end{definition}

Using the  holomorphic Lefschetz fixed-point formula 
(\cite{AS}) one can rewrite this definition as follows.

\begin{theorem}{\cite{BL1}} {\rm Let $TX \vert_{X^{g,h}}= \oplus W_{\lambda}$ and let  
$x_{\lambda}$  be the collection of Chern roots of $W_{\lambda}$. 
Let $$\Phi(g,h,\lambda,z,\tau, x)=
{{\theta({x \over 2 \pi i }+\lambda(g)-\tau \lambda(h)-z)}
 \over 
{\theta({x\over 2 \pi i }+\lambda(g)-\tau \lambda(h))}}
e^{2\pi i z \lambda(h)}.$$
Then:
\begin{equation}  \label{commute}
Ell_{orb}(X,G,z,\tau)= {1 \over {\vert G \vert}} 
\sum_{gh=hg} \prod_{\lambda(g)=\lambda(h)=0}
x_{\lambda}\prod_{\lambda} \Phi(g,h,\lambda,z,\tau,x_\lambda)[X^{g,h}]. 
\end{equation}}
\end{theorem}

The orbifold elliptic genus so defined specializes
for $q=0$ into
$$Ell_{orb} (X,G,0,y)=y^{-{\rm dim} X}\chi_{-y}(X,G)$$
where 
$$\chi_y(X,G)=\sum_{\{ g \}, X^g} y^{F(g,X^g \subset X)} 
\sum_{p,q} (-1)^q{\rm dim}H^{p,q}(X^g)^{C(g)}$$
On the other hand $\chi_y(X,G)$ is the value of the 
orbifold $E$-function 
$$E(u,v,G)=\sum_{\{g\},X^g} (uv)^{F(g,X^g \subset X)}\sum_{p,q}
{\rm dim}H^{p,q}(X^g)^{C(g)}u^pv^q$$ for $u=y,v=-1$.
In paticular $Ell_{orb} (X,G,0,1)$ coinsides with the orbifold 
Euler characteristic:
 $ e_{orb}(X,G)={1 \over {\vert G \vert}}\sum_{fg=gf} e(X^{f,g})$. 

\subsection{Discrete torsion}

\begin{definition}
Let $\alpha \in H^2(G,U(1))$, and define
\[
\delta(g,h) = \frac{\alpha(g,h)}{\alpha(h,g)}
\]
 The orbifold elliptic genus with discrete torsion $\alpha$, written
$Ell^{\alpha}_{orb} (X,G,q,y)$, is defined as
\begin{equation} \label{orbellgen}
Ell^{\alpha}_{orb} (X,G,q,y) := 
y^{-\dim  X /2} \sum_{[g],X^g}  y^{F(g,X^g\subseteq X)} \frac 1 {|C(g)|}
\sum_{h\in C(h)} \delta(g,h) L(h, V_{g,X^g\subseteq X}).
\end{equation}
As above, using the holomorphic 
Lefschetz fixed-point formula this can be rewritten as
\[
Ell^{\alpha}_{orb} (X,G;y,q) =  {1 \over {\vert G \vert}} 
\sum_{gh=hg} \delta(g,h)  \prod_{\lambda(g)=\lambda(h)=0}
x_{\lambda}\prod_{\lambda} \Phi(g,h,\lambda,z,\tau,x_\lambda)[X^{g,h}].
\]
\end{definition}

Such twisted elliptic genus has
specialization 
properties similar to the case $\alpha=0$. Using Dolbeault cohomology
corresponding to the inner local systems ${L}_{\alpha}$ 
defined by $\alpha$ 
(cf. \cite{R}) one can define twisted $E$-function:
$$E^{\alpha}(u,v,G)=\sum_{\{g\},X^g} (uv)^{F(g,X^g \subset X)}\sum_{p,q}
{\rm dim}H^{p,q}(X^g,{L}_{\alpha})^{C(g)}u^pv^q$$
which for $u=1,v=-1$ yields:
$$e^{\alpha}(X,G)=
{1 \over {\vert G \vert}}\sum_{fg=gf} \delta(f,g)e(X^{f,g})$$
The elliptic genus \ref{orbellgen} satisfies:
$$Ell^{\alpha}(0,y,G)=y^{{{\rm dim} X} \over 2}E^{\alpha}(y,-1,G)$$

We proceed to investigate the modularity properties of this twisted
orbifold elliptic genus.

\subsubsection{Jacobi forms} 
Let $\mathbb{H}$ denote the upper half plane. A \emph{weak Jacobi
form} of weight $k \in \mathbb{Z}$ and index $r \in \frac{1}{2}
 \mathbb{Z}$ is a function on $\mathbb{H} \times \mathbb{C}$ satisfying
the transformation property
\[
\phi({{a\tau +b} \over {c \tau +d}},{z \over {c \tau +d}})= (c \tau
 +d)^k e^{2 \pi i{{r cz^2} \over {c \tau +d}} }\phi (\tau, z), \; \;
 \; \; \left( \begin{matrix} a & b \\ c & d \end{matrix} \right) \in
 SL(2,\mathbb{Z})
\]
\[
\phi(\tau, z+m\tau+n)=
(-1)^{2r(m+n)} e^{-2 \pi i r (m^2 \tau +2 m z)} \phi (\tau, z),  \; \; \; \; (m,n) \in \mathbb{Z}^2
\]
and has a Fourier expansion $\sum_{l,m}c_{m,l}y^lq^m$ with nonnegative
$m$. Equivalently, we can say that a Jacobi form is an automorphic
form for the Jacobi group $\Gamma^{J} = SL(2,\mathbb{Z}) \ltimes
\mathbb{Z}^{2}$ generated by the four transformations:
$$(z,\tau) \to (z+1,\tau),
~(z,\tau) \to (z+\tau, \tau),
~(z,\tau)\to (z,\tau +1), 
~(z,\tau)\to (\frac z\tau,-\frac {1}\tau).$$

\begin{theorem}
Let $X$ be a compact complex manifold such that $K_X$ is trivial, $G$
a finite group acting effectively on $X$, and $\alpha \in
H^{2}(G,U(1))$. Let $n$ denote the order of $G$ in $\on{Aut} (H^{0}
(X, K_X))$.  Then $Ell^{\alpha}_{orb}(X,G)$ is a weak Jacobi form of
weight $0$ and index $d/2$ with respect to subgroup of the Jacobi
group $\Gamma^J$ generated by transformations
$$(z,\tau) \to (z+n,\tau),
~(z,\tau) \to (z+n\tau, \tau),
~(z,\tau)\to (z,\tau +1), 
~(z,\tau)\to (\frac z\tau,-\frac {1}\tau).$$
In particular, if the action preserves holomorphic volume then 
$Ell_{orb}(X,G)$ is a weak Jacobi form of weight $0$ and index 
$d/2$ for the full Jacobi group.

\end{theorem}

\begin{proof} 
We use the notation $Ell^{\alpha}_{orb}(X,G,z, \tau)$ rather than
$Ell^{\alpha}(X,G,q,y)$ to emphasize the dependence on $\tau$ and $z$.
It is shown in \cite{BL1} that
\begin{equation*}
\Phi(g,h,\lambda,z+1,\tau,x) =-\ee^{2  \pi \ii \lambda(h)} \cdot 
\Phi(g,h,\lambda,z,\tau,x) 
\end{equation*}
which implies that 
\[
\prod_{\lambda} \Phi(g,h,\lambda,z +n ,\tau,x_\lambda)[X^{g,h}] =
(-1)^{dn} \ee^{2 n \pi \ii \sum \lambda(h)} \prod_{\lambda}
\Phi(g,h,\lambda,z ,\tau,x_\lambda)[X^{g,h}]
\]
Now, $n \sum \lambda(h) \in \mathbb{Z}$ by the assumption that $h^n$
acts trivially on $H^{0}(X, K_X)$.  Thus 
$$Ell^{\alpha}_{orb} (X,G,
z+n, \tau) = (-1)^{dn} Ell^{\alpha}_{orb}(X,G,z,\tau).$$
The following formulas are also obtained in \cite{BL1}:
\begin{equation} \label{form1}
\Phi(g,h,\lambda,z,\tau+1,x)=\Phi(gh^{-1},h,\lambda,z,\tau,x) 
\end{equation}
\begin{equation} \label{form2}
\Phi(g,h,\lambda,z+n\tau,\tau,x)=
(-1)^n
\ee^{-2 \pi \ii nz -\pi \ii n^2 \tau}\ee^{nx+2 \pi \ii n\lambda(g)} 
\cdot 
\Phi(g,h,\lambda,z,\tau,x)
\end{equation}
\begin{equation} \label{form3}
\Phi(g,h,\lambda,\frac z \tau,\ -\frac 1 \tau,\frac x \tau) = e^{{\pi
\ii z^2 \over \tau}-\frac {zx}\tau} \cdot\Phi(h,g^{-1},\lambda,z,\tau,
x).
\end{equation}
Equation \ref{form1} and $\delta(g,h) = \delta(gh^{-1}, h)$ imply that
\[
Ell^{\alpha}_{orb}(X,G,z,\tau+1)=Ell^{\alpha}_{orb}(X,G,z,\tau).
\]
Equation \ref{form2} implies that
\[
Ell^{\alpha}_{orb}(X,G,z+n \tau,\tau)=(-1)^{dn} e^{-2 \pi \ii dn z - \pi \ii dn^2 \tau}Ell^{\alpha}_{orb}
(X,G,z,\tau)
\]
In order to see how 
\begin{equation} \label{form4}
Ell^{\alpha}_{orb}(X,G,z / \tau,-1 / \tau) = e^{\frac{\pi i d z^2}{\tau}}
Ell^{\alpha}_{orb}(X,G,z, \tau)
\end{equation}
follows, we write
\[
 \prod_{\lambda(g)=\lambda(h)=0} x_{\lambda}\prod_{\lambda}
\Phi(g,h,\lambda,z,\tau,x_\lambda) = \sum_{\bf{k}_{\lambda}}
Q(g,h,z,\tau) {\bf x_{\lambda}}^{{\bf k_{\lambda}}}
\]
where $\bf{k_{\lambda}}$ are multiindices and ${\bf x_{\lambda}}$ are
the corresponding monomials. We thus obtain,
$$
 \prod_{\lambda(g)=\lambda(h)=0} \frac{x_{\lambda}}{\tau} \prod_{\lambda}
\Phi(g,h,\lambda,\frac{z}{\tau},\frac{-1}{\tau},\frac{x_\lambda}{\tau})
= \sum_{\bf{k}_{\lambda}} (\frac{1}{\tau})^{deg(\mathbf{k}_{\lambda})}
Q(g,h,\frac{z}{\tau},\frac{-1}{\tau}) {\bf x_{\lambda}}^{{\bf
k_{\lambda}}}
$$
whereas \ref{form3} implies that 
\begin{align*}
 \prod_{\lambda(g)=\lambda(h)=0} \frac{x_{\lambda}}{\tau}
\prod_{\lambda}
\Phi(g,h,\lambda,\frac{z}{\tau},\frac{-1}{\tau},\frac{x_\lambda}{\tau})
&= e^{\frac{\pi i d z^2}{\tau}} 
\prod_{\lambda(g)=\lambda(h)=0}
\frac{x_{\lambda}}{\tau} \prod_{\lambda}
\Phi(h,g^{-1},\lambda,z,\tau, x) \\
&= \tau^{-dim(X^{g,h})} \sum_{\bf{k}_{\lambda}}
Q(h,g^{-1},z,\tau) {\bf x_{\lambda}}^{{\bf k_{\lambda}}}
\end{align*}
Thus for multiindices $\mathbf{k}_{\lambda}$ such that
$deg(\mathbf{k}_{\lambda}) = dim(X^{g,h})$, we find
\[
Q_{k}(g,h,\frac{z}{\tau}, \frac{-1}{\tau}) = Q_{k}(h,g^{-1},z,\tau)
\]
Finally, $\delta(g,h) = \delta(h,g^{-1})$ ensures that \ref{form4} holds. 
\end{proof}

\section{Discrete torsion and the chiral de Rham complex}

Let $X$ be a smooth complex algebraic variety, and $G$ a finite group
acting effectively on $X$.  In this section, we briefly review the
construction of the chiral de Rham complex of an orbifold
introduced in \cite{FS}. Another
construction of this object was obtained independently
by A. Vaintrob. 

\subsection{Vertex algebras and twisted modules} \label{VA}

In this section we will use the language of vertex superalgebras, their
modules, and twisted modules. For an introduction to vertex algebras
and their modules \cite{FLM,KAC,FBZ}, and for background on twisted
modules, see \cite{FFR,Dong,DLM,FS}.

We recall that a conformal vertex superalgebra is a
$\mathbb{Z}_+$--graded super vector space $$V = \bigoplus_{n=0}^\infty
V_n,$$ $$V_{n}=V_{n}^{\ol{0}} \oplus V_{n}^{\ol{1}}$$ together with a
vacuum vector $\vac \in V_0^{\ol{0}}$, an even translation operator $T$
of degree $1$, a conformal vector $\omega \in V_2^{\ol{0}}$ and an even
linear map

\begin{align*}
Y: V &\to \on{End} V[[z^{\pm 1}]], \\ A &\mapsto Y(A,z) = \sum_{n \in
\Z} A_{(n)} z^{-n-1}.
\end{align*}
These data must satisfy certain axioms
(see \cite{FLM,KAC,FBZ}). In what follows we will denote the
collection of such data simply by $V$, and the parity of an element $A
\in V$ homogeneous with respect to the $\mathbb{Z}/2 \mathbb{Z}$ grading by
$p(A)$.

A vector superspace $M$ is called a $V$--module if it is equipped
with an even linear map
\begin{align*}
Y^M: V &\to \on{End} M[[z^{\pm 1}]], \\
A &\mapsto Y^M(A,z) = \sum_{n \in \Z} A^M_{(n)} z^{-n-1}
\end{align*}
such that for any $v \in M$ we have $A^M_{(n)} v = 0$ for large enough
$n$. This operation must satisfy the following axioms:

\begin{itemize}

\item $Y^M(\vac,z) = \on{Id}_M$;

\item For any $v \in M$ and homogeneous $A,B \in V$ there exists an
element $$f_v \in M[[z,w]][z^{-1},w^{-1},(z-w)^{-1}]$$ such that the
formal power series
$$Y^M(A,z) Y^M(B,w) v, \quad (-1)^{p(A)p(B)}Y^{M}(B,w) Y^{M}(A,z)v ,
\qquad \on{and} $$ $$ \qquad Y_M(Y(A,z-w) B,w) v$$ are expansions of $f_v$
in $M((z))((w))$, $M((w))((z))$ and $M((w))((z-w))$, respectively.

\end{itemize}

The power series $Y^M(A,z)$ are called vertex operators. We write the
vertex operator corresponding to $\omega$ as
\[
        Y^M(\omega,z) = \sum_{n \in \mathbb{Z}} L^M_{n} z^{-n-2},
\]
where $L^M_n$ are linear operators on $V$ generating the Virasoro
algebra. Following \cite{Dong}, we call $M$ \emph{admissible} if
$L^{M}_{0}$ acts semi-simply with integral eigenvalues.

Now let $\sigma_{V}$ be a conformal automorphism of $V$, i.e., an even
automorphism of the underlying vector superspace preserving all of the
above structures (so in particular $\sigma_{V}(\omega) = \omega$). We
will assume that $\sigma_{V}$ has finite order $m>1$. A vector space
$M^\sigma$ is called a $\sigma_V$--{\em twisted} $V$--module (or
simply twisted module) if it is equipped with an even linear map
\begin{align*}
Y^{M^\sigma}: V &\to \on{End} M^\sigma[[z^{\pm \frac{1}{m}}]], \\ A
&\mapsto Y^{M^{\sigma}}(A,\zf{m}) = \sum_{n \in \frac{1}{m}\Z}
A^{M^\sigma}_{(n)} z^{-n-1}
\end{align*}
such that for any $v \in M^{\sigma}$ we have $A^{M^\sigma}_{(n)} v =
0$ for large enough $n$. Please note that we use the notation
$Y^{M^{\sigma}}(A,\zf{m})$ rather than $Y^{M^{\sigma}}(A,z)$ in the
twisted setting. This operation must satisfy the following axioms (see
\cite{FFR,Dong,DLM,LI,FS}):

\begin{itemize}

\item $Y^{M^{\sigma}}(\vac,\zf{m}) = \on{Id}_{\Mt}$;

\item For any $v \in \Mt$ and homogeneous $A,B \in V$, there exists an
element
\[
f_{v} \in \Mt [[\zf{m}, \wf{m} ]][z^{- \frac{1}{m}}, w^{-
\frac{1}{m}},(z-w)^{-1}]
\]
such that the formal power series $$Y^{\Mt}(A,\zf{m}) Y^{\Mt}(B,
\wf{m})v, \quad (-1)^{p(A)p(B)} Y^{\Mt}(B,\wf{m}) Y^{\Mt}(A, \zf{m})v,
\qquad \on{and} $$ $$ \qquad Y^{\Mt}(Y(A,z-w) B, \wf{m})v$$ are expansions
of $f_{v}$ in $\Mt((\zf{m}))((\wf{m}))$, $\Mt((\wf{m}))((\zf{m})) $
and $\Mt((\wf{m}))((z-w))$, respectively.

\item If $A \in V$ is such that $\sigma_V(A) = e^{\frac{2\pi i k}{m}}
A$, then $A^{M^\sigma}_{(n)} = 0$ unless $n \in \frac{k}{m} + \Z$.

\end{itemize}

The series $Y^{M^\sigma}(A,z)$ are called twisted vertex operators.
In particular, the Fourier coefficients of the twisted vertex operator
\[
        Y^{M^\sigma}(\omega,\zf{m}) = \sum_{n \in \mathbb{Z}}
        L^{M^\sigma}_{n} z^{-n-2},
\]
generate an action of the Virasoro algebra on $M^\sigma$. 

\subsection{The chiral de Rham complex of an orbifold}


For $g \in G$, let $X^{g}$ denote the fixed-point set of $g$, and denote by
\[
i^{g}: X^g \hookrightarrow X
\]
the inclusion map of $X^g$ into $X$.
The following results were obtained in \cite{FS}.

\begin{itemize}

\item For each $g \in G$, there exists a sheaf $\Omega^{ch,g}_X$
supported on $X^{g}$. When $g=1$, this is a sheaf of vertex superalgebras,
originally introduced in \cite{MSV}, and called the chiral de Rham
complex. We denote it simply by $\Omega^{ch}_X$. Being a sheaf of
vertex superalgebras means that for each open $U \in X$, $\Omega^{ch}_X
(U)$ is a vertex superalgebra. 

\item Let $g,h \in G$, and let $g' = h g h^{-1}$. There exist isomorphisms of sheaves
\begin{equation} \label{eqstruct}
R^{h}_{g,hgh^{-1}} : \Omega^{ch,g}_{X} \mapsto h^{*} \Omega^{ch,hgh^{-1}}_X
\end{equation}
satisfying 
\[
        R^{k}_{g',g''} \circ R^{h}_{g,g'} = R^{kh}_{g,g''}
\]
where $k \in G$ and $g''=khgh^{-1}k^{-1}$.

\item When $g \neq 1$, $\Omega^{ch,g}_X$ is a sheaf of $g$--twisted
modules, meaning that for each $g$--invariant $U$, $\Omega^{ch,g}_X
(U)$ is a $g$--twisted $\Omega^{ch}_X (U)$--module. This structure
induces a corresponding twisted module structure on
$H^{*}(X,\Omega^{ch,g}_X)$.

\item $\Omega^{ch,g}_X$ carries a bigrading by two operators $L^{g}_0,
J^{g}_0$. This bigrading induces a bigrading on
$H^{*}(X,\Omega^{ch,g}_X)$.

\item $\Omega^{ch,g}_X$ carries a differential $Q^{g}$, such that
$(Q^{g})^2 = 0$. Furthermore, there exists an inclusion of the de Rham complex of $X^{g}$
\[
i^{g}_{*} (\Omega_{dR}(X^{g}, d)) \hookrightarrow (\Omega^{ch,g}_X,Q^{g})
\]
which is a quasiisomorphism. This implies in particular that
\[
\mathbb{H}(\Omega^{ch,g}_X,Q^{g}) \cong H^{*}_{dR}(X^{g},\mathbb{C})
\]
\ref{eqstruct} implies that $C(g)$, the centralizer of $g$, acts on
$(\Omega^{ch,g}_X, Q^g)$, and therefore on its hypercohomology. This gives an isomorphism
\[
    \mathbb{H}(\Omega^{ch,g}_X,Q^{g})^{C(g)} \cong
    H^{*}_{dR}(X,\mathbb{C})^{C(g)} \cong H^{*}_{dR}(X^{g} / C(g),
    \mathbb{C})
\]
We therefore have
\begin{equation} \label{CRisom}
\bigoplus_{[g]} \mathbb{H}(\Omega^{ch,g}_X,Q^{g})^{C(g)} \cong
\bigoplus_{[g]} H^{*}_{dR}(X^g / C(g), \mathbb{C})
\end{equation}
The right-hand side is isomorphic as a
vector space to the Chen-Ruan orbifold cohomology of $[X/G]$.
Furthermore, the operators $J^g_{0}$ acting on $\Omega^{ch,g}_X$, induce a
gradation on the left which coincides with the Chen-Ruan gradation
shifted by the fermionic shift (see \cite{CR,Z}). \ref{CRisom} is therefore an
isomorphism of graded vector spaces.

\item There exists an increasing exhaustive filtration on $\Omega^{ch,g}_X$
\begin{equation} \label{filtration}
F^{0} \Omega^{ch,g}_X \subset F^{1} \Omega^{ch,g}_X \subset F^{2} \Omega^{ch,g}_X \subset \cdots.
\end{equation}
Let $\overline{\Omega}^{ch,g}_X$ denote the restriction of
$\Omega^{ch,g}_X$ to $X^g$, which inhertis a filtration from
\ref{filtration}. The bigrading operators $J^{g}_{0}, L^{g}_{0}$ are
compatible with \ref{filtration}, and so the associated graded $gr_F
(\overline{\Omega}^{ch,g}_X)$ can be described in terms of its
decomposition into eigenbundles for $J^{g}_{0}, L^{g}_{0}$. We have:

\begin{multline}    \label{Vg}
        gr_F (\overline{\Omega}^{ch,g}_X)  = \bigotimes_{k \geq 1} \left(
        \Lambda^{\bullet}_{yq^{k-1}} V^{*}_{0} \otimes
        \Lambda^{\bullet}_{y^{-1}q^{k}} V_{0} \otimes
        \on{Sym}^{\bullet}_{q^{k}} V^{*}_{0} \otimes
        \on{Sym}^{\bullet}_{q^{k}} V_{0} \otimes \right. \\
        \left.
        \bigotimes_{\lambda \ne 0} (\Lambda^{\bullet}_{yq^{k-1 +
        \lambda(g)}} V^{*}_{\lambda} \otimes
        \Lambda^{\bullet}_{y^{-1}q^{k-\lambda(g) }} V_{\lambda}
        \otimes \on{Sym}^{\bullet}_{q^{k-1+\lambda(g)}}
        V^{*}_{\lambda} \otimes
        \on{Sym}^{\bullet}_{q^{k-\lambda(g)}} V_{\lambda})
        \right).
\end{multline}
where
\[
TX \vert_{X^{g}} = \bigoplus V_{\lambda}.
\]
If we now form
\[
\mc{H}_{orb} (X,G) = \bigoplus_{[g]} H^{*} (X, \Omega^{ch,g}_{X})^{C(g)}
\]
then as shown in \cite{FS}
\[
Ell_{orb} (X,G,q,y) = \on{Supertrace}(q^{L_0} y^{J_0 - \on{dim} (X)/2}, \mc{H}_{orb} (X,G)).
\]
\end{itemize}

\subsubsection{Adding discrete torsion to the chiral de Rham complex}

In this section we show how to incorporate discrete torsion in the
above setup. Suppose that $Y$ is a $G$--manifold, and $W$ a
$G$--equivariant sheaf on $Y$. This means that for each $g \in
G$, we are given an isomorphism
\[
T_{g} : W \mapsto g^{*} W 
\]
such that 
\[
T_{g} T_{h} = T_{gh}
\]
Suppose now that $\chi: G \mapsto \mathbb{C}^{\times}$ is a character
of $G$. Then 
\[
T'_{g} = \chi(g) T_{g}
\]
is a new $G$ equivariant structure on $W$. 

We apply this observation to the sheaves $\Omega^{ch,g}_X$ and $C(g)$
rather than $G$. A class $\alpha \in H^{2} (G, U(1))$ yields characters
\[
\alpha_{g} : C(g) \mapsto U(1)
\]
defined by $\alpha_{g}(h) = \delta(g,h)$. We can now twist the
$C(g)$--equivariant structure on $\Omega^{ch,g}_X$ described above,
given by the $R^{h}_{g,g'}$, to obtain a new $C(g)$--equivariant
structure, denoted $C(g)_{\alpha}$.The following theorem is an
immediate consequence of the above discussion.

\begin{theorem}
i) 
\[
\mathbb{H}^{*} (\Omega^{ch}_X, Q_{BRST})^{C(g)_{\alpha}} \cong H^{*}_{dR} (X^{g}/C(g),
\mc{L}^{\alpha}_g)
\]
where the right-hand side is the de Rham cohomology of $[X^g/C(g)]$
with values in the orbifold local system $\mc{L}^{\alpha}_g$ described
in \cite{R}. Thus we have an isomorphism of graded vector spaces
\[
\bigoplus_{[g]} \mathbb{H}^{*} (\Omega^{ch}_X,
Q_{BRST})^{C(g)_{\alpha}} \cong H^{*}_{orb,\alpha} ([X/G], \mathbb{C})
\]
where the right-hand side is the $\alpha$--twisted Chen-Ruan cohomology of $[X/G]$ (see \cite{R}). \\
ii) Let 
\[
\mc{H}^{\alpha}_{orb}(X,G) = \bigoplus_{[g]} H^{*} (X, \Omega^{ch,g}_X)^{C(g)}_{\alpha}
\]
Then
\[
Ell^{\alpha}_{orb} (X,G,q,y) = \on{Supertrace}(q^{L_{0}} y^{J_{0} - \on{dim} (X)/2},\mc{H}^{\alpha}_{orb}(X,G))
\]
\end{theorem}

\section{Symmetric products and discrete torsion}

\subsection{The spin double cover of $S_N$}

We begin by reviewing discrete torsion for the symmetric group
following \cite{Di}.

Let $S_N$ denote the symmetric group on $N$ letters. It is
well-known that for $N\geq 4$
\begin{equation}
H^2(S_N,U(1)) \cong \Z_2.
\end{equation}
which implies that for $N \geq 4$, there is a unique non-trivial
central extension of the permutation group
\begin{equation}
1 \to \Z_2 \to \widehat S_N \to S_N \to 1.
\end{equation}

The extension $\widehat{S}_N$ can be constructed as follows. $S_N$
acts on the hyperplane in $\mathbb{R}^N$ given by 
\[
x_{1} + \cdots x_N = 0
\] 
preserving the standard inner product. This yields an embedding
\[
S_N \hookrightarrow O(N-1).
\]
Now, $O(N-1)$ has a double cover $\on{Pin}(N-1)$. Pulling back this
central extension to $S_N$ yields $\widehat{S}_N$. We call the latter
the \emph{spin double cover} of $S_N$.

In terms of generators and relations, $\widehat{S}_N$ can be described
as follows. It is generated by elements $1,z, \hat t_1, \cdots, \hat
t_{N-1}$, where $z$ is central, subject to the relations:
\begin{eqnarray}
z^2 \is 1, \nonumber \\[2mm]
\hat t_i^2 \is z,\nonumber\\[2mm]
\hat t_i \hat t_{i+1} \hat t_i \is \hat t_{i+1} \hat t_i \hat
t_{i+1},\\[2mm]
\hat t_i \hat t_j \is  z \, \hat t_j \hat t_i,\qquad j>i+1. \nonumber
\end{eqnarray}
The map $\widehat{S}_N \mapsto S_N$ amounts to sending $\hat t_i$ to
the transposition $t_i$ interchanging the $i$th and $i+1$st letters,
and sending $z$ to $1$. We can think of $z$ as being $-1$.

\subsection{Generating functions} \label{symtors}

Suppose that the elliptic genus of $X$ is given by 
\begin{equation} 
Ell(X;q,y)=
\sum_{m,\ell} c(m,\ell) q^m y^\ell 
\end{equation}
 As shown in \cite{BL1, DMVV}, the
generating function of the orbifold elliptic genera of the symmetric
products is 
\begin{equation} 
Z(p,q,y) = \sum_{N \geq 0} p^N Ell_{orb}(X^{N},S_N,q,y) =
\prod_{n>0,\,m,\ell} (1-p^nq^m y^\ell)^{-c(nm,\ell)} 
\end{equation} 
In this
section, we obtain a formula for the generating function of elliptic
genera of symmetric products with discrete torsion. Let
\begin{equation}\label{Z++}
Z^{\alpha}(p,q,y) = \sum_{N \geq 0} p^N Ell^{\alpha}_{orb}(X^{N},S_N,q,y)
\end{equation}
and let
\begin{eqnarray}
Z_{++}(p,q,y) \is  \prod_{n>0,\, m, l \geq 0} 
{\Bigl( 1 + p^{2n} q^{m - {1\over 2}}y^\ell\Bigr)^{c(n(2m-1),\ell)} \hspace{-17mm} 
\over 
\Bigl( 1 - p^{2n-1} q^{m}y^\ell\Bigr)^{c((2n-1)m,\ell)}\hspace{-17mm}} \nonu
Z_{+-}(p,q,y) \is  \prod_{n>0,\, m, l \geq 0} 
{\Bigl( 1 - p^{2n} q^{m - {1\over 2}}y^\ell\Bigr)^{c(n(2m-1),\ell)} \hspace{-17mm} 
\over 
\Bigl( 1 - p^{2n-1} q^{m}y^\ell\Bigr)^{c((2n-1)m,\ell)}\hspace{-17mm}} \nonu
Z_{-+}(p,q,y) \is  \prod_{n>0,\, m, l \geq 0} 
{\Bigl( 1 + p^{2n} q^{m}y^\ell\Bigr)^{c(2nm,\ell)} \hspace{-12mm}
\over 
\Bigl( 1 - p^{2n-1} q^{m}y^\ell\Bigr)^{c((2n-1)m,\ell)}\hspace{-17mm}} \nonu
Z_{--}(p,q,y) \is - \prod_{n>0,\, m, l \geq 0} 
{\Bigl( 1 - p^{2n} q^{m}y^\ell\Bigr)^{c(2nm,\ell)} \hspace{-12mm}
\over 
\Bigl( 1 - p^{2n-1} q^{m}y^\ell\Bigr)^{c((2n-1)m,\ell)}\hspace{-17mm}} 
\end{eqnarray}

\begin{theorem} \label{bigformula}
\begin{equation}
Z^{\alpha}(p,q,y) = {1\over 2} \left(Z_{++} + Z_{+-} + 
Z_{-+} + Z_{--}\right).
\end{equation}
\end{theorem}

We begin by recalling a variation on Lemma 4.5 from \cite{BL1} 

\begin{lemma} \label{symprod}
Let $V = V_{0} \oplus V_{1}$ be a super vector space, and $A$ and $B$
two commuting operators acting semisimply on $V$ and preserving the
parity decomposition of $V$. Assume furthermore that $B$ only has
non-negative eigenvalues in $\frac{1}{2} \mathbb{Z}$, and that the bigraded pieces
$V_{m,l} = \{ v \in V \vert Av = lv, Bv = mv \}$ are
finite-dimensional. Let $d(m,l) = sdim(V_{m,l})$, where $sdim$ denotes
superdimension.  Define the superdimension of $V$ with respect to
$A,B$ to be the series
\[
\chi(V) (y,q) = \on{Supertrace}(V,y^{A} q^{B}) = \on{tr} (V_{0},
y^{A} q^{B} ) - \on{tr} (V_{1}, y^{A} q^{B} ) = \sum_{m,l} d(m,l) q^{m} y^{l}
\]
Let $Sym^{N} V$ denote the $N$th supersymmetric product of $V$. The
operators $A$ and $B$ act on $Sym^{N} V$, and
\[
\sum_{N} p^{N} \on{Supertrace} (Sym^{N} V, y^{A} q^{B}) = \prod_{m,l}
\frac{1}{(1 - p q^{m} y^{l})^{d(m,l)}}
\]
where the right hand side is expanded in a power series in $q$ and $p$.
\end{lemma}

\medskip

Let $\Lambda^{N} V$ denote the $N$th supersymmetric wedge product of
$V$.  Since $\Lambda^{N} V$ is isomorphic to $Sym^{N} \overline{V}$,
where $\overline{V}$ denotes $V$ with its parity reversed 
(or directly from the argument in the proof of lemma 4.5 in \cite{BL1}), 
we obtain the following:

\begin{corollary} \label{wedgeprod}
Let $V$ be as in Lemma \ref{symprod}. Then
\[
\sum_{N} p^{N} \on{Supertrace} (\Lambda^{N} V, y^{A} q^{B}) = \prod_{m,l}
(1 - p q^{m} y^{l})^{d(m,l)}
\]

\end{corollary}
\bigskip

\begin{proof}(Of Theorem \ref{bigformula}) \\
Let $S_N$ denote the symmetric group on $N$ letters. We recall that
conjugacy classes in $S_N$ are parametrized by partitions of $N$. The
conjugacy class of an element $g \in S_N$ is therefore uniquely
determined by the numbers $a_j$ of $j$--cycles in the cycle
decomposition of $g$. Recall moreover that the centralizer of an
element with cycle type $[g] = (1)^{a_1} (2)^{a_2} \cdots (k)^{a_k}$
is
\[
\prod^{k}_{i=1} S_{a_i} \ltimes (\mathbb{Z} / i \mathbb{Z})^{a_i} 
\]
where the $\mathbb{Z} / i \mathbb{Z}$ act by powers of the $i$--cycles
and $S_{a_i}$ permutes the $i$--cycles among themselves. 

Let $c_j \in S_{j}$ be a $j$--cycle, and denote
$\Omega^{ch,c_j}_{X^{j}}$ simply by $\Omega^{ch,j}_{X^{j}}$. Recall
that this is a sheaf on $X^{j}$ supported on $X$ diagonally embedded,
whose fibers are twisted modules over the chiral de Rham vertex
algebra. Let $\mc{H}_{[g]} = H^{*}(X^{N}, \Omega^{ch,g}_{X^{N}})$ and
$\mc{H}_j = H^{*}(X^{j}, \Omega^{ch,j}_{X^{j}})$, viewed as a super
vector space where the parity is given by the sum of the cohomology
index and the fermionic charge grading. Furthermore, introduce the operator $D$ which acts on $\mc{H}_j$ by multiplication by $-j \on{dim}(X)/2$. 
 We have
\[
\Omega^{ch,g}_{X^{N}} = \boxtimes^{k}_{j=1}
(\Omega^{ch,j}_{X^{j}})^{\boxtimes a_j}
\]
and so by the Kunneth formula
\[
\mc{H}_{[g]} = \bigotimes^{k}_{j=1} \mc{H}_{j}^{\otimes a_j}
\]
We have 
\[
 \sum_{N} p^{N} Ell^{\alpha}_{orb} (X^{N}, S_{N},q,y)  = \sum_{N} p^N
 \sum_{[g] \in S_N} 
 \on{Supertrace} (q^{L_0} y^{J_0 + D}, \mc{H}^{C(g)_{\alpha}}_{[g]})
\]
where the subscript on $C(g)_{\alpha}$ indicates that invariants are
being taken with respect to the $\alpha$-twisted action of $C(g)$.  If
$h \in C(g)$, and $T_h$ denotes the operator of $h$ acting on
$\mc{H}_{[g]}$ untwisted by $\alpha$, then the $\alpha$--twisted
action is given by $\delta(g,h) T_{h}$.

As explained in for instance \cite{Di}, the value $\delta(g,h)$
depends on the parity of $g$ and $h$, where the latter is given by
\begin{align*}
p(g) & = \sum^{k}_{j=1} (j-1) a_j \; \; \on{mod} 2 \\
     & = \sum^{k}_{j=1, \; j \, \text{even}} a_j \; \; \on{mod} 2
\end{align*}
$C(g)$ is generated by transpositions $\tau(j)_{ab}$ interchanging two
cycles of length $j$, as well as the $j$--cycles $c_j$ in the cycle
decomposition of $g$ (we use the short-hand notation $c_j \in g$).
The result is as follows:
\begin{equation} \label{transp}
\delta(g,\tau(j)_{ab}) = (-1)^{j-1}
\end{equation}
and if $c_j \in g$, then
\begin{equation} \label{cycle}
\delta(g,c_j) = \left\{ \begin{matrix} 1, \; \; \text{if}\, j \, \text{is odd, } \\
 (-1)^{p(g)-1} \; \; \text{if} \, j \, \text{is even} \end{matrix} \right.
\end{equation}
It follows from \ref{transp} and \ref{cycle} that 
\[
\mc{H}^{C(g)_{\alpha}}_{[g]} = \bigotimes^{k}_{j = 1 \; \text{odd}}
\on{Sym}^{a_j} (\mc{H}_j^{\mathbb{Z} / j \mathbb{Z}_{\alpha}}) \otimes
\bigotimes^{k}_{j=1 \; \text{even}} \Lambda^{a_j}
(\mc{H}_j^{\mathbb{Z} / j \mathbb{Z}_{\alpha}} )
\]
The space $\mc{H}^{\mathbb{Z}/ j \mathbb{Z}_{\alpha}}$ will depend on
how $\alpha$ twists the $\mathbb{Z}/ j \mathbb{Z}$--action. For $j$
odd, there is only one possibility, and
\[
\mc{H}^{\mathbb{Z}/ j \mathbb{Z}_{\alpha}}_j = \mc{H}^{\mathbb{Z}/ j \mathbb{Z}}_j
\]
When $j$ is even, let
\begin{align*}
\mc{H}^{+}_j &= \mc{H}^{\mathbb{Z}/ j \mathbb{Z}_{\alpha}}_j, \; \; \text{when} \, [g] \, \text{is odd} \\
\mc{H}^{-}_j &= \mc{H}^{\mathbb{Z}/ j \mathbb{Z}_{\alpha}}_j, \; \; \text{when} \, [g] \, \text{is even}
\end{align*}

It was shown in \cite{BL1} that with $\mc{H} = \mc{H}^{\mathbb{Z}/ j \mathbb{Z}}_j$ or $\mc{H}^{+}_j$ 
\begin{align} \label{notors}
\on{Supertrace}(q^{L_{0}} y^{J_0 + D}, \mc{H} ) & = \frac 1j
 \sum_{r=0}^{j-1} Ell(X,q^{\frac 1j} \xi^r,y) \nonumber \\ & = \frac
 1j \sum_{m,l} c(m,l)
(\sum^{j-1}_{r=0} \xi^{mr}) q^{\frac{m}{j}} y^{l} \nonumber \\ &
 = \sum_{m,l}c(mj,l)y^lq^m.
\end{align}
where $\xi = \exp(2 \pi i/j)$. 
Similarly, using the holomorphic Lefschetz fixed-point
formula, one finds
\begin{align} \label{tors}
\on{Supertrace}(q^{L_{0}} y^{J_0 + D}, \mc{H}^{-}_j) & = \frac 1j
\sum_{r=0}^{j-1} (-1)^{r} Ell(X,q^{\frac 1j} \xi^r,y) \nonumber \\& = \frac 1j
\sum_{m,l} c(m,l)(\sum^{j-1}_{r=0} (-1)^{m} \xi^{mr}) q^{\frac{m}{j}} y^{l}
\nonumber \\ & = \sum_{m>0,l \ge 0}c((m-\frac{1}{2})j,l)y^lq^{m - \frac{1}{2}}.
\end{align}
Let
\[
\mathbb{S} = \bigotimes_{j \, \text{odd}} \on{Sym}_{p^j} \mc{H}^{\mathbb{Z}/ j \mathbb{Z}}
\]
and let
\begin{align}
\Lambda^{+} & \in \bigotimes_{j \, \text{even}} \Lambda_{p^{j}} \mc{H}^{+}_j \\
\Lambda^{-} & \in \bigotimes_{j \, \text{even}} \Lambda_{p^{j}} \mc{H}^{-}_j 
\end{align}
denote the subspaces corresponding to permutations of odd (resp. even) parity. We have that
\[
 \sum_{N} p^{N} Ell^{\alpha}_{orb} (X^{N}, S_n, q, y)  = 
 \on{Supertrace}(q^{L_0} y^{J_0 + D}, \mathbb{S} \bigotimes \Lambda^{+}) 
 +  \on{Supertrace} (q^{L_0} y^{J_0 + D}, \mathbb{S} \bigotimes \Lambda^{-}) 
\]
The result now follows from Lemma \ref{symprod}, Corollary \ref{wedgeprod}, and the observation that
\begin{align*}
\on{Supertrace}(q^{L_0} y^{J_0 + D}, \Lambda^{+}) & =\frac{1}{2} \on{Supertrace}(q^{L_0}
y^{J_0 + D}, \bigotimes_{j \, \text{even}} \Lambda_{p^{j}} \mc{H}^{+}_j) \\ & - \frac{1}{2}
\on{Supertrace}(q^{L_0} y^{J_0 + D}, \bigotimes_{j \, \text{even}}
\Lambda_{-p^{j}} \mc{H}^{+}_j)
\end{align*}
and
\begin{align*}
\on{Supertrace}(q^{L_0} y^{J_0 + D}, \Lambda^{-}) & = \frac{1}{2}
\on{Supertrace}(q^{L_0} y^{J_0 + D}, \bigotimes_{j \, \text{even}}
\Lambda_{p^{j}} \mc{H}^{-}_j) \\  & + \frac{1}{2} \on{Supertrace}(q^{L_0} y^{J_0 + D},
\bigotimes_{j \, \text{even}} \Lambda_{-p^{j}} \mc{H}^{-}_j).
\end{align*}
\end{proof}

{\it Remark.} There is an equivariant version of the 
theorem \ref{bigformula} which is also a twisted form of 
the product formula for the generating functions for the wreath 
products (conjectured in \cite{WZ} and proven in \cite{BL1} Remark 4.6. p.341).
Let $X$ and $G$ be as above and let 

\begin{equation}\label{equivar} 
Ell(X,G;q,y)=
\sum_{m,\ell} c_G(m,\ell) q^m y^\ell 
\end{equation}

The wreath product $G \wr S_N$  
(consisting of pairs $((g_1,...,g_N);\sigma), g_i \in G, \sigma \in S_N$ 
with multiplication:
$((g_1,...,g_N);\sigma_1) \cdot ((h_1,...,h_N);\sigma_2)=
((g_1 \cdot h_{\sigma_1^{-1}(1)},...,g_N \cdot h_{\sigma_1^{-1}(N)});
\sigma_1\sigma_2 )$ acts on the symmetric products $X^N$. 
The nontrivial class in $H^2(S_N,U(1))$ can be pulled
back to the class in $H^2(G \wr S_N,U(1))$ which we denote 
again as $\alpha$. 
Then the generating function 
$
Z^{\alpha}(X,G,p,q,y) = \sum_{N \geq 0} 
p^N Ell_{orb}^{\alpha}(X^{N},G \wr S_N,q,y)$ is given 
by the theorem \ref{bigformula} with the coefficients 
$c(m,l)$ in the formulas 
\ref{Z++} being replaced by $c_G(m,l)$ from \ref{equivar}.

\end{document}